\documentclass[draft]{amsart}

\usepackage{latexsym,amsmath, amssymb, amsfonts, amscd, amsthm, mathrsfs, enumerate}

\DeclareMathOperator{\aut}{Aut}

\DeclareMathOperator{\pgl}{PGL}
\newcommand{\floor}[1]{{\lfloor #1 \rfloor}}

\newtheorem{theorem}{Theorem}[section]
\newtheorem{lemma}[theorem]{Lemma}

\newtheorem{corollary}[theorem]{Corollary}

\theoremstyle{definition}
\newtheorem{remark}[theorem]{Remark}

\makeatletter
\def\emppsubsection{\@startsection{subsection}{2}{\z@}{-3.25ex plus -1ex minus -.2ex}{-1em}}

\makeatother

\newcommand \CM {{\mathcal M}}

\newcommand \CAS {{\mathcal AS}}
\newcommand \CH {{\mathcal H}}
\newcommand \PP {{\mathbb P}^1}
\newcommand \ZZ {{\mathbb Z}}
\newcommand \NN {{\mathbb N}}
\def\calc{{\mathcal C}}
\def\cald{{\mathcal D}}
\newcommand  \FF {{\mathbb F}}
\newcommand \Aut {\mathop{\rm Aut}}
\newcommand \dime {\mathop{\rm dim}}
\newcommand \Jac {\mathop{\rm Jac}}
\newcommand \Spec {\mathop{\rm Spec}}
\def\ra{\rightarrow}
\def\iso{\cong}

\newcommand{\abs}[1]{{\left|#1\right|}}

\newcommand{\st}[1]{\{#1\}}

\title{Curves of given $p$-rank with trivial automorphism group}
\author{Jeff Achter}
\address{Colorado State University, Fort Collins, CO 80523}
\email{j.achter@colostate.edu}
\author{Darren Glass} \address{Gettysburg College \\ 200 N. Washington St \\ Gettysburg, PA 17325} \email{dglass@gettysburg.edu}
\author{Rachel Pries$^{*}$} \address{Colorado State University, Fort Collins, CO 80523}
\email{pries@math.colostate.edu}


\thanks{* The third author was
partially supported by NSF grant DMS-07-01303.}
\date{}

\begin{document}

\begin{abstract}
Let $k$ be an algebraically closed field of characteristic $p >0$.
Suppose $g \geq 3$ and $0 \leq f \leq g$. We prove there is a smooth
projective $k$-curve  of genus $g$ and $p$-rank $f$ with no non-trivial
automorphisms. In addition, we prove there is a smooth projective
hyperelliptic $k$-curve  of genus $g$ and $p$-rank $f$ whose
only non-trivial automorphism is the hyperelliptic involution. The
proof involves computations about the dimension of the moduli space
of (hyperelliptic) $k$-curves of genus $g$ and $p$-rank $f$ with extra automorphisms.
\end{abstract}

\maketitle

\section{Introduction}

Let $k$ be an algebraically closed field of characteristic $p >0$.  If
$g \geq 3$, there exist a $k$-curve $C$ of genus $g$ with
$\Aut(C)=\{1\}$ and a hyperelliptic $k$-curve $D$ of genus $g$ with
$\Aut(D) \simeq \ZZ/2$ (see, e.g., \cite{Popp} and \cite{fischer56},
respectively).  In this paper, we extend these results to curves with
given genus and $p$-rank.

If $C$ is a smooth projective $k$-curve of genus $g$ with Jacobian $\Jac(C)$,
the {\it $p$-rank} of $C$ is the integer $f_C$ such that the cardinality
of $\Jac(C)[p](k)$ is $p^{f_C}$.
It is known that $0 \leq f_C \leq g$.  We prove the following:

\begin{corollary} \label{C1}
Suppose $g \geq 3$ and $0 \leq f \leq g$.
\begin{enumerate}[(i)]
\item There exists a smooth projective $k$-curve $C$ of genus $g$ and $p$-rank $f$
with $\Aut(C)=\{1\}$.
\item There exists a smooth projective hyperelliptic $k$-curve $D$ of genus $g$ and $p$-rank $f$
with $\Aut(D) \simeq \ZZ/2$.
\end{enumerate}
\end{corollary}

More generally, we consider the moduli space $\CM_g$ of curves of
genus $g$ over $k$.
The $p$-rank induces a stratification $\CM_{g,f}$ of $\CM_g$
so that the geometric points of $\CM_{g,f}$
parametrize $k$-curves of genus $g$ and $p$-rank at most $f$.
Similarly, we consider the $p$-rank stratification $\CH_{g,f}$ of
the moduli space $\CH_g$ of hyperelliptic $k$-curves of genus $g$.
Our main results (Theorems \ref{T1} and \ref{T2}) state that, for every geometric generic point $\eta$ of $\CM_{g,f}$
(resp.\ $\CH_{g,f}$), the corresponding curve $\calc_\eta$ satisfies $\Aut(\calc_\eta)=\{1\}$
(resp.\ $\Aut(\cald_\eta) \simeq \ZZ/2$).

For the proof of the first result, we consider the locus
$\CM_g^{\ell}$ of $\CM_g$ parametrizing $k$-curves of genus $g$ which have an
automorphism of order $\ell$. Results from \cite{FVdG:complete} and
\cite{Popp} allow us to compare the dimensions of $\CM_{g,f}$ and
$\CM_{g}^{\ell}$. The most difficult case, when $\ell=p$,
involves wildly ramified covers and deformation results from
\cite{BM}. For the proof of the second result, we compare the dimensions of $\CH_{g,f}$ and
$\CH_{g}^{\ell}$ using \cite{GP} and \cite{GV}. When $p=2$, this
relies on \cite{PZ}. The hardest case for hyperelliptic curves is when $p \geq 3$, $f=0$, and
$\ell=4$ and we use a degeneration argument to finish this case.

The statements and proofs of our main results would be simpler if
more were known about the geometry of $\CM_{g,f}$ and $\CH_{g,f}$.
For example, one could reduce to the case $f=0$ if one knew that
each irreducible component of $\CM_{g,f}$ contained a component of
$\CM_{g,0}$. Even the number of irreducible
components of $\CM_{g,f}$ (or $\CH_{g,f}$) is known only in
special cases.

We also sketch a second proof of the main results that uses degeneration to the
boundaries of $\CM_{g,f}$ and $\CH_{g,f}$, see Remark \ref{Radv1}.

\begin{remark}
  There is no information in Corollary \ref{C1} about the field of
  definition of the curves. In the literature, there are several
  results about curves with trivial automorphism group which are
  defined over finite fields. In \cite{Po} and \cite{Po2}, the author
  constructs an $\FF_p$-curve $C_0$ of genus $g$ with
  $\Aut_{\bar\FF_p}(C_0)=\{1\}$ and a hyperelliptic $\FF_p$-curve
  $D_0$ of genus $g$ with $\Aut_{\bar\FF_p}(D_0) \simeq \ZZ/2$.
  However, the $p$-ranks of $C_0$ and $D_0$ are not considered.

For $p=2$ and $0 \leq f \leq g$, the author of \cite{Z} constructs a
hyperelliptic $\FF_2$-curve $D_0$ of genus $g$ and $p$-rank $f$ with
$\Aut_{\bar\FF_p}(D_0) \simeq \ZZ/2$. The analogous
question for odd characteristic appears to be open.  Furthermore, for
all $p$ it
seems to be an open question whether there exists an $\FF_p$-curve
$C_0$ of genus $g$ and $p$-rank $f$ with
$\Aut_{\bar\FF_p}(C_0)=\{1\}$ \cite[Question 1]{Z}.
\end{remark}

\subsection{Notation and background}

All objects are defined over an algebrai\-cally closed field $k$ of
characteristic $p >0$. Let $\CM_g$ be the moduli space of smooth
projective connected curves of genus $g$, with tautological curve
$\calc_g \ra \CM_g$. Let $\CH_g$ be the moduli space of smooth
projective connected hyperelliptic curves of genus $g$, with
tautological curve $\cald_g \ra \CH_g$.

If $C$ is a $k$-curve of genus $g$, the $p$-rank of $C$ is the
number $f\in\st{0, \ldots, g}$ such that $\Jac(C)[p](k) \iso
(\ZZ/p)^f$.  The $p$-rank is a discrete
invariant which is lower semicontinuous in families. It induces a
stratification of $\CM_g$ by closed reduced subspaces $\CM_{g,f}$
which parametrize curves of genus $g$ with $p$-rank at most
$f$. Similarly, let $\CH_{g,f}\subset \CH_g$ be the locus of
hyperelliptic curves of genus $g$ with $p$-rank at most $f$.

Recall that $\dime(\CM_g)=3g-3$ and $\dime(\CH_g)=2g-1$.
Every irreducible component of $\CM_{g,f}$ has dimension $2g-3+f$ by \cite[Thm.\ 2.3]{FVdG:complete}.
Every irreducible component of $\CH_{g,f}$ has dimension $g-1+f$ by \cite[Thm.\ 1]{GP} when $p \geq 3$ and
by \cite[Cor.\ 1.3]{PZ} when $p=2$.
In other words, the locus of curves of genus $g$ and $p$-rank $f$ has
pure codimension $g-f$ in $\CM_g$ and in $\CH_g$.

Every irreducible component of $\CM_{g,f}$ (and $\CH_{g,f}$) has a geometric generic point $\eta$.
Let $\calc_\eta$ (resp.\ $\cald_\eta$) denote the curve corresponding to the point $\eta$.

Let $\ell$ be prime.  Let $\CM_g^\ell \subset \CM_g$ denote the locus
of curves which admit an automorphism of order $\ell$ (after pullback by a finite cover of the base).
The locus $\CM_g^\ell$ is closed in $\CM_g$.
If $D$ is a hyperelliptic curve, let $\iota$ denote the unique
hyperelliptic involution of $D$.  Then $\iota$ is in the center of $\Aut(D)$.
Let $\CH_g^\ell \subset \CH_g$ denote the locus of hyperelliptic
curves which admit a non-hyperelliptic automorphism of order $\ell$.
Let $\CH_{g}^{4,\iota}$ denote
the locus of hyperelliptic curves which admit an automorphism $\sigma$
of order four
such that $\sigma^2=\iota$.

An Artin-Schreier curve is a curve that admits a structure as
$\ZZ/p$-cover of the projective line.
Let $\CAS_{g} \subset \CM_g$ denote the locus of Artin-Schreier curves of genus $g$ and
let $\CAS_{g,f}$ denote its $p$-rank strata.

Unless stated otherwise, we assume $g \geq 3$ and $0 \le f \le g$.
\section{The case of $\CM_g$}

\subsection{A dimension result}

Suppose $\Theta$ is an irreducible component of $\CM_{g}^{\ell}$ with generic point $\xi$.
Let $Y$ be the quotient of $\calc_{\xi}$ by a group of order $\ell$.
Let $g_Y$ and $f_Y$ be respectively the genus and $p$-rank of $Y$.
Consider the $\ZZ/\ell$-cover $\phi:\calc_{\xi} \to Y$.
Let $B\subset Y$ be the branch locus of $\phi$.
If $\ell=p$, let $j_b$ be the jump in the lower ramification
filtration of $\phi$ at a branch point $b \in B$ \cite[IV]{Se:lf}.

\begin{lemma} \label{Lemmadim}
\begin{enumerate}[(i)]
\item
If $\ell \not = p$, then $\dime(\Theta) \leq 2(g-g_Y)/(\ell-1)+f_Y-1$;
\item
If $\ell=p$, then $\dime(\Theta) \leq 2(g-g_Y)/(\ell-1)+f_Y-1-\sum_{b \in B} \lfloor j_b/p \rfloor$.
\end{enumerate}
\end{lemma}

\begin{proof}
Let $\phi:\calc_\xi \to Y$ be as above, with branch locus $B\subset
Y$.  Because $g \ge 3$, if $g_Y = 1$ then $\abs B > 0$.
Let $\CM_{g_Y,f_Y,|B|}$ be the moduli space of curves of genus $g_Y$
and $p$-rank at most $f_Y$ with $|B|$ marked points.  Then
$\dime(\CM_{g_Y,f_Y,|B|})=2g_Y-3+f_Y+|B|$ if $g_Y \geq 1$.
Also $\dime(\CM_{0,0,|B|})=|B|-3$ if $\abs B \ge 3$.

\begin{enumerate}[(i)]
\item
Since $\phi:\calc_\xi \ra Y$ is tamely ramified, the curve
$\calc_\xi$ is determined by the quotient curve $Y$, the branch
locus $B$, and ramification data which is discrete. Therefore,
$\dim (\Theta) \le \dim (\CM_{g_Y,f_Y,\abs B})$ if $g_Y \geq 1$. and
$\dim(\Theta) \leq |B|-3$ if $g_Y =0$. By the Riemann-Hurwitz
formula, $2g-2 = \ell(2g_Y-2) + \abs B (\ell-1)$.  One can deduce
that $\abs B = 2(g-\ell g_Y)/(\ell-1) +2$ and the desired result
follows.

\item
By the Riemann-Hurwitz formula for wildly ramified covers \cite[IV, Prop.\ 4]{Se:lf},
$$2g-2=p(2g_Y-2)+\sum_{b \in B}(j_b+1)(p-1).$$

For $b \in B$, let $\hat{\phi}_b:\hat{\calc}_z \to \hat{Y}_b$ be the germ of the cover $\phi$
at the ramification point $z$ above $b$.
By \cite[p.229]{BM},
the dimension of the moduli space of covers $\hat{\phi}_b$ with ramification break $j_b$
is $d_b=j_b-\lfloor j_b/p \rfloor$.
The local/global principle of formal patching (found, for example, in \cite[Prop.\ 5.1.3]{BM}) implies
$\dime(\Theta) \le \dime(\CM_{g_Y, f_Y, |B|})+\sum_{b \in B}d_b$.
Since $|B|+\sum_{b\in B}j_b= 2(g-pg_Y)/(p-1)+2$, this simplifies to
$$\dime(\Theta)\le 2(g-g_Y)/(p-1)+f_Y-1 -\sum_{b \in B}\lfloor j_b/p \rfloor.$$
\end{enumerate}
\end{proof}

\subsection{No automorphism of order $p$}

\begin{lemma} \label{Lmgfp}
Suppose $\Gamma$ is a component of $\CM_{g,f}$ with geometric generic point $\eta$.
Then $\calc_{\eta}$ does not have an automorphism of order $p$.
\end{lemma}

\begin{proof}
The strategy of the proof is to show that $\dime(\Gamma \cap
\CM_g^p) < \dim(\Gamma)$.  Recall that $\dim(\Gamma) = 2g-3+f$ by
\cite[Thm.\ 2.3]{FVdG:complete}.

Let $\Theta$ be an irreducible component of $\Gamma\cap \CM_g^p$,
with geometric generic point $\xi$.
Consider the resulting cover
$\phi:\calc_{\xi} \to Y$, which is either \'etale or wildly
ramified. Let $g_Y$ and $f_Y$ be respectively the genus and $p$-rank
of $Y$.

Suppose first that $g_Y=0$.  In other words, $\xi \in \CAS_{g,f}$ and
$\calc_{\xi}$ is an Artin-Schreier curve.  By \cite[Lemma 2.6]{PZ},
$g=d(p-1)/2$ for some $d \in \NN$.  If $p=2$, then
$\dime(\CAS_{g,f})=g-1+f$ \cite[Cor.\ 1.3]{PZ}. If $p \geq 3$, then
$\dime(\CAS_{g,f}) \leq d-1$ by \cite[Thm.\ 1.1]{PZ}. In either case,
$\dime(\Theta) \le \dime(\CAS_{g,f}) < \dime(\Gamma)$ since $g \geq 3$.

Now suppose that $g_Y \geq 1$.
If $p \geq 3$, Lemma \ref{Lemmadim}(ii) implies that $\dime(\Theta) \leq g-g_Y+f_Y-1 < 2g-3+f$.

If $p=2$ and if $g_Y \ge 1$, let $\abs B$ be the number of branch points of $\phi$.
By the Deuring-Shafarevich formula \cite[Cor.\ 1.8]{Crew}, $f-1=2(f_Y-1)+|B|$.
Lemma \ref{Lemmadim}(ii) implies that $\dime(\Theta)\le 2g-2g_Y +(f-1-|B|)/2
-\sum_{b \in B}\lfloor j_b/2 \rfloor$.
In particular, $\dime(\Theta) < 2g-2g_Y+f/2$.
So $\dime(\Theta) < 2g-3+f$ if $g_Y \geq 2$.

Suppose $p=2$ and $g_Y=1$.
The hypothesis $g \geq 3$ implies that $\phi$ is ramified.
So $|B| \geq 1$ and $j_b \geq 1$ for $b \in B$.
Then $\dime(\Theta) < 2g-3+f/2$.

Thus $\dime(\Theta)< \dim(\Gamma)$ in all cases. This
inequality implies that $\eta \not \in \CM_g^{p}$ and
$\Aut(\calc_{\eta})$ does not contain an automorphism of order $p$.
\end{proof}

\subsection{The main result for $\CM_{g,f}$}

\begin{theorem}  \label{T1}
Suppose $g \geq 3$ and $0 \leq f \leq g$. Suppose $\eta$ is the
geometric generic point of an irreducible component $\Gamma$ of
$\CM_{g,f}$. Then $\Aut(\calc_{\eta})=\{1\}$.
\end{theorem}

\begin{proof}
By Lemma \ref{Lmgfp}, $\Aut(\calc_{\eta})$ contains no automorphism
of order $p$.  Let $\ell \not =p$ be prime.
Consider an irreducible component $\Theta \subset \Gamma \cap \CM_{g}^\ell$.
The result follows in any case where  $\dim(\Theta) < \dim(\Gamma)=2g-3+f$.

Let $\xi$ be the geometric generic point of $\Theta$.
Let $Y$ be the quotient of $\calc_\xi$ by a group of order $\ell$.
Let $g_Y$ and $f_Y$ be the genus and $p$-rank of $Y$.

If $\ell\ge 3$, then Lemma \ref{Lemmadim}(i) implies $\dime(\Theta)
\le g-g_Y+f_Y-1$.  Thus $\dime(\Theta) < 2g-3+f$ and $\calc_\eta$ has
no automorphism of order $\ell \geq 3$.

Suppose $\ell=2$. If $g_Y=0$, then $\calc_\eta$ is hyperelliptic
and in particular $\dime(\Theta) \leq \dime(\CH_{g,f})=g-1+f <
2g-3+f$. If $g_Y \geq 1$, then $\dime(\Theta) \le 2g-2g_Y+f_Y-1$
which is less than $2g-3+f$ except when $g_Y=1$ and $f=f_Y \leq
1$.

For the final case, when $\ell =2$, $g_Y = 1$, and $f = f_Y$, Lemma
\ref{Lemmadim} alone does not suffice to prove the claim.  Let
$\CM_g^{2,Y}$ be the moduli space of curves of genus $g$ which are
$\ZZ/2$-covers of $Y$.  It is the geometric fiber over the moduli point of $Y$ of a map from a proper, irreducible Hurwitz space to $\CM_1$ (see, e.g., \cite[Cor.\
6.12]{bertinromagny07}).  Therefore, $\CM_g^{2,Y}$ is irreducible.
Now $\xi \in
\CM_g^{2,Y}\cap \Gamma$.  The strategy is to show that there exists $s \in
\CM_g^{2,Y}$ such that $f_s > f_Y$.  From this, it follows that
$\CM_g^{2,Y} \cap \CM_{g,f_Y}$ is a closed subset of
$\CM_g^{2,Y}$ of positive codimension.  Then $\Theta$ is a closed subset of $\Gamma$ of positive codimension, and the proof is complete.

To construct $s$,
consider a $\ZZ/2$-cover $\psi_1:Y \to \PP$.
If $g$ is odd (resp.\ even), let $\psi_2:X \to \PP$
be a $\ZZ/2$-cover so that $X$ has genus $(g-1)/2$ (resp.\ $g/2$)
and so that the branch locus of $\psi_2$ contains exactly $2$ (resp.\ $3$)
of the branch points of $\psi_1$.
Since only $2$ (resp.\ $3$) of the branch points of $\psi_2$ are specified,
one can suppose $X$ is ordinary.
Consider the fiber product $\psi:W \to \PP$ of $\psi_1$ and $\psi_2$.
Following the construction of \cite[Prop.\ 3]{GP},
$W$ has genus $g$ and $p$-rank at least $g/2$.
Since $W$ is a $\ZZ/2$-cover of $Y$, it corresponds to a point $s \in
\CM_g^{2,Y}$ with $p$-rank at least $f_Y+1$.
\end{proof}

Here is the proof of part (i) of Corollary \ref{C1}:

\begin{corollary}
\label{C1parti}
Suppose $g \geq 3$ and $0 \leq f \leq g$.
There exists a smooth projective $k$-curve $C$ of genus $g$ and $p$-rank $f$
with $\Aut(C)=\{1\}$.
\end{corollary}

\begin{proof}
Let $\Gamma$ be an irreducible component of $\CM_{g,f}$, with
geometric generic point $\eta$.  Let $\Gamma' \subset \Gamma$ be the open, dense subset
parametrizing curves with $p$-rank exactly $f$ \cite[Thm.\
2.3]{FVdG:complete}.  By Theorem \ref{T1}, $\aut(\calc_{\eta}) =
1$. The sheaf $\underline{\aut}(\calc)$ is constructible on
$\Gamma'$, but there are only finitely many possibilities for the
automorphism group of a curve of genus $g$. Therefore, there is a
nonempty open subspace $U\subset\Gamma'$ such that, for each $s\in
U(k)$, $\calc_s$ has $p$-rank $f$ and $\aut(\calc_s) = 1$.
\end{proof}

\begin{corollary} \label{C2}
Let $g \geq 3$ and $0 \leq f \leq g$. There exists a principally
polarized abelian variety $(A, \lambda)$ over $k$ of dimension $g$
and $p$-rank $f$ with $\Aut(A, \lambda)=\{\pm 1\}$.
\end{corollary}

\begin{proof}
Let $A$ be the Jacobian of the curve given in Corollary
\ref{C1parti}.
The desired properties then follow from Torelli's theorem
\cite[Thm.\ 12.1]{Mi}.
\end{proof}

\section{The case of $\CH_g$}

Recall that $g\geq 3$ and $0 \leq f \leq g$.

\subsection{When $p=2$}

\begin{lemma} \label{Lp=2}
Let $p=2$ and suppose $\eta$ is the geometric generic point of a
component $\Gamma$ of $\CH_{g,f}$. Then $\aut(\cald_{\eta}) \simeq
\ZZ/2$.
\end{lemma}

\begin{proof}
The automorphism group of a hyperelliptic curve always
contains a (central) copy of $\ZZ/2$.  Let $U\subset \Gamma$ be the
subset parametrizing curves with automorphism group $\ZZ/2$.  As in
the proof of Corollary \ref{C1parti}, $U$ is open; it suffices to show
that $U$ is nonempty.

By \cite[Cor.\ 1.3]{PZ}, $\CH_{g,0}$ is irreducible of
dimension $g-1$ when $p=2$.  For $g \geq 3$, there exists a
hyperelliptic curve $D_0$ with $p$-rank $0$ and $\Aut(D_0) \simeq
\ZZ/2$ \cite[Thm.\ 3]{Z}. The component $\Gamma$ contains $\CH_{g,0}$ by
\cite[Cor.\ 4.6]{PZ}. Then $U$ is non-empty since $U\cap
\CH_{g,0}$ is nonempty.
\end{proof}

\subsection{No automorphism of order $p$}\label{S:wild}

Suppose $p \geq 3$.

\begin{lemma} \label{Lwild}
If $p |(2g+2)$ or $p|(2g+1)$, then $\dime \CH_g^p= \lfloor (2g+2)/p
\rfloor -2$.  Otherwise, $\CH_{g}^p$ is empty.
\end{lemma}

\begin{proof}
  Suppose $s \in \CH_g^p(k)$.  There exists
  $\sigma \in \aut(\cald_s)$ of order $p$.  Since $\iota$ and $\sigma$
  commute, $\sigma$ descends to an automorphism of $\cald_s/\langle \iota
  \rangle \simeq \PP$. Let $Z$ be the projective line $\cald_s/\langle
  \sigma, \iota \rangle$.  Then $\cald_s \to Z$ is the fiber product of the
  hyperelliptic cover $\phi:\cald_s/\langle \sigma \rangle \to Z$ and the
  $\ZZ/p$-cover $\psi: \cald_s/\langle \iota \rangle \to Z$.

Since $\cald_s/\langle \iota \rangle$ has genus zero, the cover $\psi$ is
  ramified only at one point $b$ and the jump $j_b$ in the lower ramification filtration
equals 1.  After changing coordinates on
  $\cald_s/\langle \iota \rangle$ and $Z$, the cover $\psi$ is isomorphic to
  $c^p-c=x$.

If $\phi$ is not branched at $\infty$ then each branch
  point of $\phi$ lifts to $p$ branch points of the cover $\cald_s \to
  \cald_s/\langle \iota \rangle$, and the branch locus of $\phi$ consists of
  $(2g+2)/p$ points.  On the other hand, if $\phi$ is branched at
  $\infty$ then the branch locus of $\phi$ consists of $(2g+1)/p$
  points. Therefore, if $\CH_g^p(k)$ is nonempty, then either
  $p|(2g+1)$ or $p|(2g+2)$.

Moreover, any branch locus of size
  $\floor{(2g+2)/p}$ uniquely determines such a cover $\phi$.
 A point $s \in \CH_g^p$ is determined by the
  branch locus of $\phi$ up to the action of affine linear
  transformations on $Z$. Thus $\dim(\CH_g^p) = \lfloor(2g+2)/p\rfloor
   -2$.
\end{proof}

\begin{lemma} \label{noautp}
Let $\eta$ be the geometric generic point of a component of
$\CH_{g,f}$. Then $\Aut(\cald_{\eta})$ contains no automorphism of
order $p$.
\end{lemma}

\begin{proof}
By Lemma \ref{Lwild}, $\CH_g^p$ is either empty or of dimension $\lfloor (2g+2)/p \rfloor -2$.
If $g \ge 3$, then
$\dim(\CH_g^p)<g-1+f = \dim(\CH_{g,f})$.
Thus $\cald_{\eta}$ does not have an automorphism of order $p$.
\end{proof}

\subsection{Extra automorphisms of order two and four}

Suppose $p \ge 3$.  In this section, we
show that the geometric generic point of any component of $\CH_{g,f}$
parametrizes a curve with no
extra automorphism of order two or four.
The proof relies on degeneration and requires an analysis of curves
of genus $2$ and $p$-rank $0$.

\begin{lemma} \label{Lg2f0}
Suppose $p \geq 3$ and $g=2$. If $\eta$ is a geometric generic
point of $\CH_{2,0}$, then $\aut(\cald_{\eta}) \simeq \ZZ/2$.
\end{lemma}

\begin{proof}
By \cite[p.130]{IKO}, $\Aut(\cald_\eta)/\langle \iota \rangle
\simeq G$ where $G$ is one of the following groups: $\{1\}$,
$\ZZ/5$, $\ZZ/2$, $S_3$, $\ZZ/2\oplus \ZZ/2$, $D_{12}$, $S_4$, or
$\pgl_2(\ZZ/5)$. Let $T^G \subset \CH_{2,0}$ be the sublocus
parametrizing hyperelliptic curves $D$ with $\Aut(D)/\langle \iota
\rangle \simeq G$.  Since every component of $\CH_{2,0}$ has
dimension one, it suffices to show that each $T^G$ is zero-dimensional.

If $G=\ZZ/5$ and $s \in T^G(k)$,  then the Jacobian
of $\cald_s$ has an action by $\ZZ/5$, and thus must be one of the two
abelian surfaces with complex multiplication by $\ZZ[\zeta_5]$.
Therefore, there exist at most two hyperelliptic
curves $D$ of genus $2$ and $p$-rank $0$ with $\Aut(D)/\langle
\iota \rangle \simeq \ZZ/5$.

Now let $G$ be any non-trivial group from the list other than $\ZZ/5$.  A curve
of genus two and $p$-rank zero is necessarily supersingular, and
any supersingular hyperelliptic curve $D$ of genus two with $\Aut(D)/\langle \iota \rangle \simeq G$
is superspecial by \cite[Prop.\ 1.3]{IKO}.
Since there are only finitely many superspecial abelian surfaces,
 $T^G$ is a proper closed
subset of $\CH_{2,0}$ for each $G \not = \{1\}$ on the list. Thus
$\aut(\cald_{\eta}) \simeq \ZZ/2$.
\end{proof}

\begin{lemma} \label{Lell=2}
Suppose $p \ge 3$ and  $g \geq 3$.
\begin{enumerate}[(i)]
\item Then $\CH_{g}^2$ is irreducible with dimension $g$;
\item there exists $s \in \CH_{g}^2(k)$ such that $\cald_s$ has $p$-rank at least $2$;
\item and $\dime(\CH_{g,0} \cap \CH_{g}^2) < g-1$.
\end{enumerate}
\end{lemma}

\begin{proof}
Suppose $s \in \CH_g^2(k)$.  There is a
Klein-four cover $\phi:\cald_s \to \PP_k$ such that $\phi$ is the fiber
product of two hyperelliptic covers $\psi_i:C_i \to \PP_k$ \cite[Lemma
3]{GP}.

If $g$ is even, then one can assume that $C_1$ and $C_2$ both have
genus $g/2$ and that the branch loci of $\psi_1$ and $\psi_2$ differ
in a single point.  If $g$ is odd, then one can assume that $C_1$
has genus $(g+1)/2$, $C_2$ has genus $(g-1)/2$, and the branch locus
of $\psi_2$ is contained in the branch locus of $\psi_1$
\cite[Prop.\ 3]{GP}. In both cases, the third $\ZZ/2$-subquotient of
$\cald_s$ has genus zero.  In particular, if $f_s$ denotes
the $p$-rank of $\cald_s$ then $f_s = f_{C_1}+ f_{C_2}$ \cite[Cor.\
2]{GP}.

\begin{enumerate}[(i)]
\item This is found in \cite[Cor.\ 1]{GP}.
\item One can choose $\psi_1$ so that $C_1$ is ordinary.  Then $f_s
  \ge \lceil\frac g2 \rceil \ge 2$.

\item  Suppose $s \in \CH_{g,0}(k)$, so that $f_s = f_{C_1} = f_{C_2}
  = 0$.
If $g$ is even, then the parameter space for choices of $\psi_1$ has
dimension $\dime (\CH_{g/2,0}) = g/2-1$.  For fixed $\psi_1$,
the parameter space for choices of $\psi_2$ has dimension at
most $1$.  Similarly, if $g$ is odd, the parameter space for choices of $\psi_1$
has dimension $\dime(\CH_{(g+1)/2,0}) = (g-1)/2$.  For fixed $\psi_1$,
there are at most finitely many possibilities for $\psi_2$.
In either case $\dime(\CH_{g,0} \cap \CH_g^2) \le \floor{g/2} < g-1$.
\end{enumerate}
\end{proof}

\begin{lemma} \label{Lell4}
Suppose $p \geq 3$ and $g \geq 3$.  Then $\CH_{g}^{4, \iota}$ is
irreducible with dimension $g-1$ and its geometric generic point
parametrizes a curve with positive $p$-rank.
\end{lemma}

\begin{proof}
Suppose $s \in \CH_g^{4,\iota}(k)$.
Let $\sigma$ be an automorphism of $\cald_s$ of order $4$ such that
$\sigma^2=\iota$. Consider the $\ZZ/4$-cover $\cald_s \stackrel{\alpha}{\to}
\PP_x \stackrel{\beta}{\to} \PP_z$. Then $\beta$ is branched at two points
and ramified at two points. Without loss of generality, one can
suppose these are $0_x$ and $\infty_x$ on $\PP_x$ and $0_z$ and
$\infty_z$ on $\PP_z$. This implies that the action of $\sigma$ on
$\PP_x$ is given by $\sigma(x)=-x$.

The inertia groups of $\beta \circ \alpha$ above $0$ and $\infty$ are
subgroups of $\langle \sigma \rangle \simeq \ZZ/4$ which are not
contained in $\langle \sigma^2 \rangle$.  Thus they each have order
4 and $\alpha$ will be branched over $0_x$ and $\infty_x$. The other $2g$
branch points of $\alpha$ form orbits under the action of $\sigma$ and
one can denote them by $\{\pm \lambda_1, \ldots, \pm \lambda_g\}$.
Without loss of generality, one can suppose $\lambda_g=1$ and
$\beta(\lambda_g)=1$ and therefore $\cald_s$ has an affine equation of the
form $y^2=x(x^2-1)\prod_{i=1}^{g-1}(x^2-\lambda_i^2)$.

Let $S=\PP-\{0,1,\infty\}$.  Let $\Delta \subset S^{g-1}$ be the strong diagonal consisting of all
$(g-1)$-tuples $(x_1, \ldots, x_{g-1})$ so that $x_i=x_j$ for some $i \not = j$.
Let $\Delta' \subset S^{g-1}$ consist of all $(g-1)$-tuples $(x_1, \ldots, x_{g-1})$
so that $x_i =- x_j$ for some $i \not = j$.
There is a surjective morphism $\omega:(\PP-\{0,1,\infty\})^{g-1}-(\Delta  \cup \Delta') \to \CH_{g}^{4, \iota}$,
where
$\omega$ sends $(\lambda_1, \ldots, \lambda_{g-1})$ to the isomorphism class of the curve with affine equation
$y^2=x(x^2-1)\prod_{i=1}^{g-1}(x^2-\lambda_i^2)$.
Thus $\CH_g^{4, \iota}$ is irreducible.

There are only finitely many fractional linear transformations fixing
the set $\{\pm \lambda_1, \ldots, \pm \lambda_{g-1}, \pm 1, 0, \infty\}$.
Thus $\omega$ is finite-to-one and $\dime(\CH_{g}^{4, \iota})=g-1$.

Suppose $g \geq 3$, and let $\eta$ be the geometric generic point of
$\CH_g^{4,\iota}$.  To finish the proof, it suffices to show that the
$p$-rank of $\cald_\eta$ is positive.
Let $T=\Spec(k[[t]])$ and let $T'=\Spec(k((t)))$.
Consider the image of the $T'$-point $(t \lambda_1, t\lambda_2,
\lambda_3, \ldots, \lambda_{g-1})$ under $\omega$. This gives a
$T'$-point of $\CH_g^{4, \iota} \subset \CH_{g}$.
The moduli space $\overline{\CH}_g$ of stable
hyperelliptic curves is proper,
so the $T'$-point of $\CH_g$  gives rise to a $T$-point of $\overline{\CH}_g$.
The special fiber of this $T$-point corresponds to a stable curve $Y$.
The stable curve $Y$ has two components $Y_1$ and $Y_2$
intersecting in an ordinary double point.
Here $Y_1$ has genus $2$ and affine equation $y_1^2=x(x^2-\lambda_1^2)(x^2-\lambda_2^2)$, while
$Y_2$ has genus $g-2$ and affine equation
$y_2^2=\prod_{i=3}^{g-1}(x^2-\lambda_i^2)$.

The moduli point $s\in \overline{\CH}_g(k)$ of $Y$ is in the closure
of $\CH_g^{4,\iota}$.  The automorphism $\sigma$ extends to $Y$, and
stabilizes each of the two components $Y_1$ and $Y_2$.  Therefore, the
moduli point of $Y_1$ lies in $\CH_2^{4,\iota}$.  There is a
one-parameter family of such curves $Y_1$ since one can vary the
choice of $\lambda_2$.  By Lemma \ref{Lg2f0}, one can suppose that
$f_{Y_1} \not =0$.  Now $f_Y =f_{Y_1}+f_{Y_2}$ by \cite[Ex.\
9.2.8]{blr}.  Thus $f_{Y} \not =0$.  Since the $p$-rank can only
decrease under specialization, and since $s$ is in the closure of
$\eta$, the $p$-rank of $\cald_\eta$ is non-zero as well.
\end{proof}

\subsection{Main Result for $\CH_{g,f}$}

\begin{theorem}
\label{T2}
Suppose $g \geq 3$ and $0 \leq f \leq g$.
If $\eta$ is the geometric generic point of an irreducible component
of $\CH_{g,f}$,
then $\Aut(\cald_{\eta}) \simeq \ZZ/2$.
\end{theorem}

\begin{proof}
Let $\Gamma$ be the irreducible component of $\CH_{g,f}$ whose
geometric generic point is $\eta$. Suppose $\sigma \in
\Aut(\cald_{\eta})$ has order $\ell$ with $\sigma \not \in \langle
\iota \rangle$. Then $p \geq 3$ by Lemma \ref{Lp=2}. Without loss
of generality, one can suppose that either $\ell$ is prime or
$\ell=4$ with $\sigma^2=\iota$.

If $\ell=4$ and $\sigma^2=\iota$, then $\CH_g^{4, \iota}$ is
irreducible with dimension $g-1$ by Lemma \ref{Lell4}. This is
strictly less than $\dime(\Gamma)$ unless $f=0$. If $f=0$, the
two dimensions are equal but the geometric generic point of $\CH_g^{4, \iota}$
corresponds to a curve of non-zero $p$-rank by Lemma \ref{Lell4}. Thus $\cald_{\eta}$
has no automorphism $\sigma$ of order $4$ with $\sigma^2=\iota$.

If $\ell$ is prime, one can suppose that $\ell \not =p$ by Lemma
\ref{noautp}. In \cite[p.10]{GV}, the authors use an argument similar
to the proof of Lemma \ref{Lwild} to show that $\CH_g^\ell$ is empty
unless $\ell \mid (2g+2-i)$ for some $i \in \{0,1,2\}$;
and if $\CH_g^\ell$ is non-empty then its dimension is
$d_{g,\ell}=-1 + (2g+2-i)/\ell$.
If $d_{g,\ell} < \dime(\Gamma)=g+f-1$
then $\cald_{\eta}$ cannot have an automorphism of order $\ell$.
This inequality is always satisfied when $\ell \geq 3$ since $g \geq 3$.

Suppose $\ell=2$.  Then
$d_{g,\ell} < \dime(\Gamma)$ unless $f \leq 1$. If $f = 1$ then
the two dimensions are equal.
By Lemma \ref{Lell=2}, $\CH_g^2$ is irreducible and contains the moduli point of a curve with
$p$-rank at least two. Therefore, the component $\Gamma$ of
$\CH_{g,1}$ is not the same as the unique irreducible component of $\CH_g^2$.

Finally, suppose $\ell=2$ and $f=0$. By Lemma \ref{Lell=2}(iii),
$\dime(\Gamma \cap \CH_{g,0})< g-1$. Thus $\eta\not\in \CH_g^2$, and
$\aut(\cald_{\eta}) \simeq \ZZ/2$.
\end{proof}

Part (ii) of Theorem \ref{C1} now follows:

\begin{corollary}
\label{C1partii} Suppose $g \ge 3$ and $0 \le f \le g$.  There
exists a smooth projective hyperelliptic $k$-curve $D$ of genus $g$
and $p$-rank $f$ with $\aut(D) \simeq \ZZ/2$.
\end{corollary}

\begin{proof}
The result
follows from Theorem \ref{T2}, using the same argument that was used to
deduce Corollary \ref{C1parti} from Theorem \ref{T1}.
\end{proof}

\begin{remark} \label{Radv1}
The proof of the last statement of Lemma \ref{Lell4} uses the
intersection of $\overline{\CH}_g^{4,\iota}$ with
the boundary component $\Delta_2$ of $\overline{\CH}_g$.
More generally, one can give a different proof of the main results of this paper using induction.
Here are the main steps of the inductive proof.
If $g \geq 3$ and if $1 \leq i \leq g/2$,
one can show that the closure of every component of $\CM_{g,f}$ in $\overline{\CM}_g$
intersects the boundary component $\Delta_i$ by \cite[p.80]{Di},
\cite{L:complete}.
Points of $\Delta_i$ correspond to singular curves $Y$ that have two components $Y_1$ and $Y_2$
of genera $i$ and $g-i$ intersecting in an ordinary double point.
Using a dimension argument, one can show that $Y_1$ and $Y_2$ are generically smooth and that their $p$-ranks
$f_1$ and $f_2$ add up to $f$.
If the generic point of a component of $\CM_{g,f}$ parametrizes a
curve with a nontrivial automorphism, another dimension argument shows that
this automorphism stabilizes each of $Y_1$ and $Y_2$.
This would imply that the generic point of a component of $\CM_{g-i,f_2}$
parametrizes a curve with nontrivial automorphism group, which would
contradict the inductive hypothesis.

An analogous proof works for $\CH_{g,f}$ when $p \geq 3$ using \cite{FVdG:complete}.

One can also use monodromy techniques to prove Corollary \ref{C2}, see \cite[App.\ 4.4]{achterpries07b}.
\end{remark}

\bibliographystyle{hamsplain}
\bibliography{paperAUT}

\end{document}